\documentclass[matzei]{svjour}
\begin{document}
\newcommand{\und}{\underline}
\newcommand{\arr}{\Rightarrow}
\newcommand{\li}{\tilde}
\newcommand{\q}{^{[q]}}
\newcommand{\m}{{\bf m}}
\newcommand{\aaa}{{\bf a}}
\newcommand{\bbb}{{\bf b}}
\newcommand{\ccc}{{\bf c}}
\newcommand{\param}{\underline{x}}
\newcommand{\tpar}{\underline{x}^{[t]}}
\newcommand{\tparq}{\underline{x}^{[tq]}}
\newcommand{\bs}{\boldsymbol}
\newcommand{\cor}{^{<q>}}
\newcommand{\lar}{\longrightarrow}
\title{Tight closure and linkage classes in Gorenstein rings}
\author{Adela Vraciu}
\thanks{I thank Mel Hochster and Craig Huneke for their support and
encouragement and for many helpful discussions.}
\institute{ Department of Mathematics, University of Kansas, Lawrence, KS
  66045}
\date{Received: date / Revised version: date}
\maketitle
\begin{abstract}
{We study the relationship between the tight closure of an
ideal and the sum of all ideals in its linkage class}
\end{abstract}
\maketitle

\section{Introduction}
\label{intro}
Tight closure and linkage have been developed as independent branches of
commutative algebra. Tight closure was introduced by Hochster and
Huneke in \cite{HH1} as a tool for studying ideals in rings of
positive characteristic, while the theory of linkage (liaison) has its
roots in the study of curves in three dimensional projective space
(\cite{Ap}, \cite{Ga}, \cite{Rao}, etc.). The algebraic foundations of linkage were
established in \cite{PS}, \cite{AN}, \cite{HU}, et cetera. 

The main result of this paper (Thm. ~\ref{main}) establishes a connection between these two
theories. In the process of proving Thm.~\ref{main} we establish a
relationship between ideals $J$ in the linkage class of an ideal $I$
and ideals in the linkage
class of an $\m $-primary ideal $(I, x_1, \ldots, x_{d-g})$, where
$d=\mathrm{dim}(R)$, $g=\mathrm{ht}(I)$, and the choice of $x_1,
\ldots, x_{d-g}$ depends on $J$
(Prop.~\ref{linkage}). This result might be of independent interest.
We also develop a theory of {\it corner powers} of unmixed ideals, which
are obtained as direct links of Frobenius powers. We explore some
of their properties in Section 3, and we use them as a tool
in the proof of Thm.~\ref{main}.

The setting is that of a Gorenstein ring of positive characteristic, where one can
use the information available about the tight closure of ideals of
finite projective dimension to  relate the tight closure of an unmixed
ideal to ideals in its linkage class.

The investigation carried out in this paper is motivated by the
following question: If $R$ is a characteristic $p$ ring with test
ideal $\tau$, and $I$ is an arbitrary ideal, what is the relationship
between $I:\tau$ and $I^*$? The definition of the test ideal
(see~\ref{test ideal}) implies
that $I:\tau \supset I^*$, and in the particular case of a Gorenstein
characteristic $p$ ring $R$ and an ideal $I$ of finite projective
dimension
we actually have equality: $I^*=I:\tau $ (see~\ref{essential}).

This equality is far from true in general. We propose to seek 
ideals that multiply $I:\tau $ into $I^*$.
Note that in the case when $I:\tau $ is multiplied into $I^*$ by an
$\m$-primary ideal, it follows that $I$ admits a test exponent, and
therefore the
tight closure of $I$ commutes with localization (see \cite{HHexp} for
details about test exponents and localization of tight closure).

The main result of this paper (Thm.~\ref{main}) states the following:
  if $R$
Gorenstein and $I$ is unmixed, then 
$(\tilde{I})^*(I:\tau )\subset I^*$, where $\tilde{I}$ denotes the sum
  of all the ideals in the linkage class of $I$. As an application, it follows that
  tight closure commutes with localization for any ideal $I$ for which $\li
  I$ is $\m $-primary (Cor.~\ref{commutes}).

We do not expect that $\tilde{I}$ is the largest ideal with
this property, but this result is interesting in light of the
unexpected relationship between tight closure and linkage.
There are several interesting consequences pertaining to properties of
the linkage class of certain ideals. For instance,
 if $I$ is an unmixed tightly closed ideal
containing the test ideal, then $I$ is maximal in its linkage class,
in the sense of containing every ideal in its linkage class
(see Cor.~\ref{max}). This provides a class of examples addressing
the question raised in \cite{PU}:  for which ideals $I$ is every ideal
in the linkage class of $I$ contained in $I$?

\section{Preliminaries}
\label{sec:1}

In this paper, $(R, m)$ denotes a Gorenstein local ring of
characteristic $p>0$ and  $q=p^e$ denotes a power of $p$.
By {\it parameter ideal} we mean an ideal
generated by part of a system of parameters.

We recall the relevant definitions:
\begin{definition}
Let $I$ be an ideal of $R$. For every $q=p^e$, $I\q:=(i^q \, | \, i \in I)$ is called the
{\it Frobenius q{\rm th} power}
     of $I$. 
$R^0$ denotes the set of elements in $R$ which are not in any minimal
prime.

An element $x \in R$ is in the {\it tight closure} of $I$ if there
exists $c \in R^0$ such that $cx^q \in I\q$ for all $q=p^e$. We write
$x \in I^*$.
\end{definition}

\begin{definition}\label{test ideal}
The {\it test ideal} of $R$ is the ideal $\tau $ generated by all
elements $c \in R^0$ such that for every ideal $I$ and every $x \in I^*$
we have $ cx \in I$.
\end{definition}

\begin{definition}
Let $(R, \m)$ be a Noetherian local ring. An ideal
$I$ of height $g$ is called {\it unmixed} if all the associated prime ideals of $I$
have height $g$.
\end{definition}

\begin{definition}
An ideal $\aaa $ is called {\it Gorenstein} if the ring $R/\aaa $ is
Gorenstein.
\end{definition}
 Note that if $R$ is a Gorenstein ring, then any parameter
ideal is a Gorenstein ideal. If $R$ is Gorenstein and $R/\aaa $ has finite projective
dimension over $R$, $R/\aaa $ is Gorenstein if and only if
its minimal free resolution over $R$ is self-dual.

The following property of unmixed ideals is well-known:
\begin{note}\label{note}
Let $(R, \m )$ be a Noetherian local ring, let $I$ be an unmixed ideal
of height $g$, and let $\aaa \subset I$ be a Gorenstein ideal of
height $g$. Then
$
\aaa :(\aaa :I)=I.
$
If the assumption that $I$ is unmixed is removed, we have $\aaa :(\aaa
:I)=I^{unm}$, where $I^{unm}$ denotes the intersection of the primary
components in a primary decomposition of $I$ associated to primes of height $g$.
\end{note}
In this paper, linkage will mean linkage by Gorenstein ideals of
finite projective dimension:
\begin{definition}
Let $I$ and $J$ be unmixed ideals of height $g$.

 We say that $I$ and
$J$ are {\it directly linked} if there is a Gorenstein ideal of finite
projective dimension $\aaa \subset I \cap J$ of the same height $g$
such that $I=\aaa :J$ (and $J=\aaa :I$).

We say that $I$ and $J$ are {\it
linked} (in $n$ steps), or that $J$ is in the {\it linkage class} of $I$, if there is a sequence of ideals $I_1=I, I_2,
\ldots, I_{n+1}=J$ such that $I_i$ and $I_{i+1}$ are directly linked
for all $i=1, \ldots, n$.
We denote the sum of all
ideals in the linkage class of $I$ by $\tilde I$.
\end{definition}

\section{Corner powers}
\label{sec:2}

The corner powers of an unmixed ideal $I$, denoted $I \cor$,  are a 
tool used to study the relationship between the tight closure of an
ideal and ideals in its linkage class; they
 are obtained as links of Frobenius powers. In this section
 we define the corner powers of unmixed ideals and study some of their properties.
\begin{definition}
Let $(R, \m)$ be a characteristic $p$ Gorenstein ring, and let $q$ be
a power of $p$.
If $I\subset R$ is an unmixed ideal of height $g$, let $\aaa \subset I$ be a
Gorenstein ideal of finite projective dimension and of height $g$ (for
example a parameter ideal) and let $J:=\aaa :I$. We define
the $q${\it th corner power of $I$} to be the ideal $I\cor :=\aaa^{[q]}:J\q$.
\end{definition}

\begin{note} Let $R, I, \aaa , J$ be as above.
If moreover $I$ is $\m $-primary, we can also write $J\q = \aaa \q :
I\cor $ by Note.~\ref{note}, since $J\q $ is unmixed and $\aaa \q $ is
Gorenstein. However, this equality is not true in general for non $\m
$-primary ideals, because
even though $J$ is unmixed, $J\q $ might not be unmixed. We can say in
general that $\aaa \q :I\cor = (J\q )^{{unm}}$.

Also note that if $I$ is $\m $-primary, we have $I^{<pq>}=(I\cor
)^{<p>}$, since $J\q =\aaa \q :I\cor \Rightarrow (I\cor )^{<p>} =
(\aaa \q)^{[p]} : (J\q)^{[p]}=\aaa ^{[pq]}:J^{[pq]} = I^{<pq>}.$ We do
not know if this equality is true for unmixed ideals which are not $\m $-primary.
\end{note}

Note that a priori the definition of $I\cor $ depends on the choice
of  $\aaa $. In order to show that $I\cor $ is well-defined
we need the following preliminary results:
\begin{lemma}\label{mapping cone} Let $R$ be a Gorenstein
local ring
  and let 
 $\aaa \subset \bbb $ be Gorenstein ideals of the
  same height $g$ and of finite projective dimension.

 The natural map $R/\aaa \lar R/\bbb $ extends to a map of complexes
 $\psi .$
 from a minimal free resolution of $R/\aaa $ to a minimal 
 free resolution of $R/\bbb$. Let $\delta \in R$ be such that the last map
$\psi _g :R
 \lar R$ is multiplication by $\delta $.

 Then $\aaa :\bbb =(\aaa, \delta )$, and 
 $\bbb   =\aaa :\delta $.
\end{lemma}
\begin{proof}
  Prop.2.6 in \cite{PS} shows that the mapping cone of the map $\psi
\check .$ is the projective free resolution of $R/(\aaa :\bbb )$.
Since the first map in the mapping cone complex is given by $(a_1, \ldots, a_g, \delta
)$, where $\aaa =(a_1, \ldots, a_g)$, it follows that $\aaa :\bbb
=(\aaa , \delta )$, and therefore $\aaa :\delta =\bbb $ by Note.~\ref{note}.
\end{proof}
\begin{lemma}\label{unm}
Let $(R, \m )$ be a Gorenstein ring, $I$ an unmixed ideal of height
$g$, and let $\aaa \subset \bbb \subset I$ be Gorenstein ideals of height
$g$ and finite projective dimension. Let $\delta $ be as in Lemma
~\ref{mapping cone}. Then the ideal $(\aaa , \delta I)$ is unmixed.
\end{lemma}
\begin{proof}
We have the following short exact sequence:
$$
0 \lar \frac{(\aaa , \delta )}{(\aaa , \delta I)} \lar \frac{R}{(\aaa ,
  \delta I)} \lar \frac{R}{(\aaa, \delta )}\lar 0.
$$
The first term in the short exact sequence is isomorphic to
$$
\frac{R}{(\aaa , \delta I):\delta}.
$$
If $u \in (\aaa :\delta I):\delta $, we can write $\delta (u-i) \in
\aaa $ for some $i \in I$, hence $ u-i \in \aaa :\delta =\bbb \subset
I$.
This shows that the first term in the short exact sequence is
isomorphic to
$R/I$, hence unmixed, while the last term is $R/(\aaa :\bbb )$, which is also
unmixed. Therefore the middle term, $R/(\aaa + \delta I)$ is unmixed,
which finishes the proof of the lemma.
\end{proof}

\begin{proposition} Let $(R, \m)$ be a Gorenstein ring of characteristic $p$.
The corner powers of an unmixed ideal $I$ are well-defined.
\end{proposition}
\begin{proof}
        Let $\aaa , \bbb \subset I$ be Gorenstein ideal of finite projective dimensions of height
$g$. Without loss of generality, we can assume $\aaa \subset
\bbb
$, since otherwise we may replace $\aaa $ by an $\aaa'\subset \aaa
\cap \bbb$.
 Let $\delta $ be as in Lemma~\ref{mapping cone}, and let 
\begin{equation}\label{def}
J_1=\aaa :I,  \ \ \ \ \ \ \  J_2=\bbb :I=(\aaa
:\delta): I=\aaa :\delta I.
\end{equation}
Note that we can also write
\begin{equation}\label{two}
I=\bbb :J_2=\aaa :\delta J_2
\end{equation}
Note that $\bbb \q $ (respectively $\aaa \q$) is again a Gorenstein ideal, and a
free resolution of $R/\bbb ^q$ (respectively $R/\aaa \q$) is obtained from
a free resolution of $R/\bbb $ (respectively $R/\aaa $)
by raising all the entries of the matrices appearing in the free resolution
to the $q$th power (by Thm. 1.7 in \cite{PS2}). Lemma~\ref{mapping
  cone} applied to $\aaa \q \subset \bbb \q $ shows that $\bbb \q
=\aaa \q :\delta ^q$.

We need to show that for all $q=p^e$, we have $\aaa \q :J_1\q =\bbb
\q:J_2\q$.
Note that $\bbb\q :J_2\q=\aaa \q: \delta ^q J_2\q$, and therefore it
suffices to show that $(\aaa \q, \delta ^q J_2\q)=(\aaa \q
,J_1\q)$. This is obviously true, since equation~\ref{two} combined with
Note.~\ref{note} shows
$(\aaa, \delta J_2)=(\aaa,
J_1)=\aaa :I$ (note that $(\aaa, J_1)= J_1$ and $J_2$ are unmixed by
construction, and $(\aaa, \delta J_2)$ is unmixed by Lemma~\ref{unm}).
\end{proof}
We would like to see how  the corner powers $I\cor $ are related to
the Frobenius powers $I\q$:
\begin{proposition}\label{cor}
Let $(R, \m)$ be a Gorenstein characteristic $p$ ring and let $I$ be an
unmixed ideal. For all $q=p^e$ we have $I\q \subset I\cor$ and if $I$
has finite projective dimension we have equality.
If
$R$ is a hypersurface the equality holds for all $q$ if and only if
$I$ has finite projective dimension.
\end{proposition}
\begin{proof}
Choose $\aaa \subset I$ a parameter ideal and
let $J=\aaa :I=(f_1, \ldots, f_n)$.
To show that $I\q \subset I\cor$ note that $I\q =(\aaa :J)\q \subset
\aaa\q :J\q$, because $I\q J\q = (\aaa :J)\q J\q =((\aaa :J) J )\q \subset \aaa
\q$.

Consider the short exact sequence
$$
0 \lar \frac{R}{I}\buildrel (f_1, \ldots, f_n)\over\lar \bigoplus _n
\frac{R}{\aaa}\lar N \lar 0$$
where $N$ is the cokernel of the map given by $(f_1, \ldots, f_n)$. Recall that 
$R^{1/q}$ denotes the $R$-algebra obtained by adjoining all $q$th
roots of elements in $R$, and it can be identified with $R$ via the
action of the Frobenius map $x \rightarrow x^q$.
Tensoring this short exact sequence with $R^{1/q}$ yields
$$
0 \lar \mathrm{Tor}_1^R(N, R^{1/q}) \lar \frac{R^{1/q}}{IR^{1/q }}\buildrel
  (f_1, \ldots, f_n)\over\lar \bigoplus _n \frac{R^{1/q}}{\aaa R^{1/q}}\lar
  N\otimes R^{1/q}\lar 0$$
and therefore
$\mathrm{Tor}_1^R(N, R^{1/q})$ is isomorphic to 
$$\frac{\aaa R^{1/q}:_{R^{1/q}} (f_1, \ldots, f_n)}{IR^{1/q}}, $$
which, upon identifying $R^{1/q}$ with $R$, can be identified with
$$\frac{\aaa \q :(f_1^q, \ldots, f_n^q)}{ I\q} =\frac{I\cor}{I\q}.$$

If $I$ has finite projective dimension, it follows that $N$ also has
finite projective dimension, and therefore Tor$_1^R(N, R^{1/q})=0$
(see Thm. 1.7 in \cite{PS2}),
and the short exact sequence shows that $I\q =I\cor $.

If $R$ is a hypersurface and $I\q =I\cor $ for all $q$, it follows
that Tor$_1^R(N, R^{1/q})=0$ for all $q$, which implies that $N$ has
finite projective dimension by \cite{He}, and therefore $I$ has
finite projective dimension.
\end{proof}

The following property of corner powers is essential for the purpose
of this paper:
\begin{theorem}\label{essential}
Let $(R, \m )$ be a Gorenstein ring of characteristic $p$, with test
ideal $\tau $. 

a). If $\aaa $ is an ideal of finite projective dimension,
then $\aaa ^*=\aaa :\tau $.

b). If $I$ is an unmixed ideal, then for all $q=p^e$ we have
$$
I\cor :\tau \supset (I:\tau )\q.
$$
More precisely:
$$
I:\tau =\lbrace x \in R \, |\,
 cx^q \in I^{<q>} \ for \ some \ c\in R^0 \ and \   all\  q=p^e
\rbrace 
$$
$$=\lbrace x \in R \, |\, cx^q \in I\cor \ for\  all \ c \in \tau \ and
\ all\ q=p^e 
\rbrace
$$
\end{theorem}
\begin{proof}
The statement in part a). for the case when $\aaa $ is a
parameter ideal is Cor.4.2(2) in \cite{Hu}. The general case will
follow from part b).

To prove part b)., choose $\aaa \subset I$ a parameter ideal of the same
height as $I$ and let $c \in \tau $.  We have:
$cx^q \in I^{<q>}\iff cx^q \in \aaa \q :J\q  \iff cx^qJ\q \subset \aaa
\q .
$
This holds for all $q$ if and only if $xJ\subset \aaa ^*=\aaa :\tau  \iff x \in
\aaa :\tau J=I:\tau $.

For the general case of part a)., use the fact that $\aaa \q =\aaa \cor
$ for any $\aaa $ of finite projective dimension (see~\ref{cor}), and
therefore for any $u \in \aaa :\tau $ we have $\tau u^q \subset \aaa \cor
=\aaa \q$, hence $u \in \aaa ^*$.
\end{proof}

\begin{proposition} \label{HK} Let $(R, m)$ be a Gorenstein local ring.
If $J$ is an $m$-primary ideal, then $l(R/J\cor )$  behaves like a Hilbert-Kunz function, that is, there is a real constant $c_J>0$ such that
$$
l\left(\frac{R}{J\cor}\right)=c_Jq^d+ {\mathcal O} (q^{d-1}),
$$
where $d={\mathrm dim}\  R.$
More precisely, we have 
$$
l(R/J\cor)=l(R/{\aaa}^{[q]})-l(R/I\q ),
$$
where $\aaa \subset J$ is a Gorenstein ideal of finite projective dimension, and $I=\aaa :J$.
\end{proposition}
\begin{proof} We have  $J\cor ={\aaa}^{[q]}:I\q$.
Note that $R/{\aaa}^{[q]}$ maps onto $R/I\q $.

Since $R$ is Gorenstein and $\aaa$ is a Gorenstein ideal of finite projective dimension, $R/\aaa \q $ is a Gorenstein ring; hence we have 
$$
\frac{R}{{\aaa}^{[q]}}=E_{R/{\aaa}^{[q]}}(k) \ \arr \frac{J\cor }{{\aaa}^{[q]}} =\mathrm{ Ann}_{R/{\aaa}^{[q]}} \left( \frac{I\q }{{\aaa}^{[q]}}\right)  =E_{R/I\q }(k),
$$
and therefore
$$
 l\left( \frac{J\cor }{{\aaa}^{[q]}} \right) =l\left( \frac{R }{I\q } \right) .$$
This can be rewritten as
$$
l\left( \frac{R }{{\aaa}^{[q]}} \right )-l\left( \frac{R}{J\cor } \right) =
l\left( \frac{R }{I\q } \right) ,
$$
and therefore $l(R/J\cor) $ can be written as a difference of
Hilbert-Kunz functions (see \cite{Mo}).
\end{proof}

\begin{proposition} \label{higher}
Let $(R,m)$ be a Gorenstein local ring with test ideal $\tau $. Assume that $R$ is excellent and analytically irreducible.

a). If $I\supset \tau$, then for all $q$,  $I\cor \supset \tau$. 

b). If $I$ is an ideal strictly contained in $\tau $, then there is a
$q_0$ such that $I\cor \subset m^{[q/q_0]}$ for all $q\ge q_0$.

\end{proposition}
\begin{proof} a). Let $\aaa \subset I $ be a Gorenstein ideal of finite projective dimension and let
$J=\aaa :I\subset \aaa :\tau$. Using the fact that
$\aaa ^*=\aaa : \tau $, and therefore $(\aaa :\tau )\q \subset \aaa \q
:\tau $, we have
 $$I\cor = \aaa \q :J\q  \supset \aaa \q :(\aaa
:\tau )\q \supset \aaa \q :(\aaa \q:\tau )\supset \tau.$$

b). Let $\aaa \subset I$ be a Gorenstein ideal of finite projective dimension, and let $J=\aaa :I$. The assumptions imply that $\aaa :\tau=\aaa^*$ is strictly contained in $J$.
Let $f\in J\backslash \aaa^*$. Then the Prop. 2.4 in \cite{Ab} shows
    that there is a $q_0$ such that 
$$
I\cor ={\aaa}^{[q]}:J\q \subset {\aaa}^{[q]}:f^q
\subset m^{[q/q_0]}
$$
\end{proof}

\begin{example}
It is not always true that $I\cor \subset I$. Let

${R=k[[X,Y,Z]]/(X^3+Y^3+Z^3)},$
$p=\mathrm{ char} (k)=2$,  $I=(x^2, y^2, z^2)$ and $ \aaa =(x^2, y^2)$. Then $$ J=(x^2, y^2):I=(x^2, y^2):z^2=(x^2, y^2, z).$$
Then $I^{<2>}=(x^4, y^4):z^2=(x^4, y^4, xyz)$, and we have $xyz \in
I^{<2>}\backslash I$.
\end{example}

The usefulness of corner powers is further illustrated in the
following proposition (compare to exercise 2.8 in \cite{Hu2} and Prop.
3.3(d) in \cite{HHexp}):
\begin{proposition}\label{decr}
Let $R$ be a Gorenstein ring of characteristic $p$ with test ideal
$\tau $ and let $I$ be an arbitrary ideal. For all $q=p^e$, let
$I_q:=\lbrace u \in R | u^q \subset \tau I\q :\tau \rbrace$.
Then for all $q$
we have $I_{pq}\subset I_{q}$, and therefore $I^*$ can be written as a
nonincreasing intersection of the ideals $I_q$.
\end{proposition}
\begin{proof}
First assume that $I$ is $\m $-primary, so that if $\aaa \subset I$ is
an $\m $-primary parameter ideal and $J=\aaa :I$ we have $I\q =\aaa
\q :J \cor $.

Let $u \in I_{pq}$, so $\tau u^{pq } \in \tau I^{[pq]}$.
Since $J^{<pq>}=\aaa ^{[pq]}:I^{[pq]}$, and $J^{<pq>}:\tau = \aaa
^{[pq]}:\tau I^{[pq]}$, this implies that
$\tau u^{pq} (J^{<pq>}:\tau ) \subset \aaa ^{[pq]}$. Using
Thm.~\ref{essential}(b) applied to $J\cor $, we get
$$
u^{pq}(J\cor :\tau )^{[p]}\subset u^{pq}(J^{<pq>}:\tau )\subset \aaa ^{pq}:\tau =(\aaa ^{[pq]})^*,
$$
which implies that
$$
u^q (J\cor :\tau ) \subset (\aaa \q )^*=\aaa \q :\tau ,
$$
and therefore $\tau u^q \in \aaa \q :(J\cor :\tau ) = \aaa \q + \tau
I\q $ (the last equality follows because $\aaa \q : (\aaa \q +\tau
I\q)= (\aaa \q :I\q ):\tau =J\cor :\tau $, and $\aaa \q +\tau I\q $ is
an unmixed ideal - being $\m $-primary -  so we can use Note.~\ref{note}).
Since this holds for any parameter ideal $\aaa \subset I$, we get
$\tau u^q
\in \tau I\q $ by Krull's intersection theorem.

If $I$ is not $\m $-primary and $u^{pq}\in \tau I^{[pq]}:\tau$, we have
$u^{pq}\in \tau (I+\m^t)^{[pq]}:\tau $ for all $t >0$, which according to
the $\m $-primary case implies \newline  $u^q \in \tau (I+\m^t)\q :\tau $ for all
$t$, hence intersecting over all $t>0$ yields $u^q \in \tau I\q :\tau
$.

In order to justify the last sentence of the proposition, note that 
$\cap_q I_q \subset I^*$ by the definition of tight closure, and vice
versa $I^*\subset I_q$ for all $q$, since $u \in I^*\Rightarrow u^q
\in (I\q )^* \subset \tau I\q :\tau $ (the last implication follows by
Thm.3.1 in \cite{V2}).
\end{proof}

\begin{corollary}\label{testexp}
Let $(R, \m)$ be a Gorenstein ring with test ideal $\tau $ and let $I$
be an ideal such that $(I:\tau )/I^*$ has finite length. Then $I$
admits a test exponent, i.e. there exists a $q_0=p^{e_0}$ such that if
$\tau x^{q_0} \in \tau I^{[q_0]}$ with $x\in R$, then $x\in I^*$.
In particular, tight closure commutes with localization for $I$.
\end{corollary}
\begin{proof}
According to Prop.~\ref{decr}, the assumption that $(I:\tau )/I^*$
has finite length implies that $I^*=I_{q_0}$ for some $q_0$, since
$I^*$ can be written as a nonincreasing intersection of the ideals
$I_q$, which are contained in $I_1=\tau I:\tau \subset I:\tau $.
\end{proof}

\begin{note}
The notion of test exponent used here is a modification of the notion
introduced in \cite{HHexp}, which basically states that $\tau x^{q_0} \in
I^{[q_0]} \Rightarrow x \in I^*$. The condition used here is weaker,
but still sufficient to guarantee that tight closure commutes with
localization in the context of Cor.~\ref{testexp}, using the fact that
in a Gorenstein local ring the test
ideal commutes with localization (Thm.4.1 in \cite{Sm4}) and it is a
strong test ideal (Thm.3.1 in \cite{V2}). The
proof follows along the lines of Prop.2.3 in \cite{HHexp}: if
$W\subset R$ is a multiplicative system and $u/1 \in (I_W)^*$, then we
have $\tau u^{q_0}/1 \in \tau (I_W)^{[q_0]}$, and we can choose $f\in W$
such that $f \tau u^{q_0} \in \tau I^{[q_0]}$, so $\tau (fu)^{q_0} \in
\tau I^{[q_0]}$. But then $fu \in I^*$, and so $u \in (I^*)_W$.
\end{note}
\section{Main result}
\label{sec:3}

The following theorem is the main result of this paper:

\begin{theorem}\label{main}
Let $(R, \m)$ be a Gorenstein local ring with test ideal $\tau $.
Let $I$ be an unmixed ideal, and let denote $\tilde{I}$ the sum 
of all ideals in the linkage class of $I$.

Then we have
$$
(\li{I})^* (I:\tau )\subset      I^*.
$$
\end{theorem}

 The outline of the proof is as follows: we prove the result for the
 case when $I$ is $\m $-primary, then we reduce to the $\m $-primary
 case
via adjoining extra elements to all the ideals involved. Our proof in the $\m $-primary
 case rests on the ability to exploit the relationships between corner
 powers and Frobenius powers of the ideal $I$, which in turn rests on
 the fact that for every $q$, $I\q $ is still an unmixed ideal (this
 is obvious if $I$ is $\m $-primary, but rarely true otherwise).

 In order to complete the announced reduction to the $\m $-primary
 case, we establish a relationship between the linkage class of an
 ideal $I$ of height $g$ and the linkage class of an $\m $-primary
 ideal containing $I$ (see Prop.~\ref{linkage}).

 To this end we need the following preliminary result:
\begin{lemma}\label{case 1}
Let $(R, \m )$ be a Gorenstein ring, $I$ an unmixed ideal of height
$g$, and let $\aaa \subset \bbb \subset I$ be Gorenstein ideals of height
$g$ and finite projective dimension.

Then we have
$$
\aaa :(\bbb :I)=\aaa +I(\aaa :\bbb).
$$
\end{lemma}
\begin{proof}
We need to check that the ideal on the right hand side is unmixed and
that the equality holds after taking duals into $\aaa $.
Indeed,
$$
\aaa:[\aaa : (\bbb :I)]=\bbb :I,
$$
and
$$
\aaa:[\aaa +I(\aaa :\bbb)]=[\aaa: (\aaa :\bbb )]:I =\bbb :I,
$$
because $\aaa \subset \bbb $.
The fact that $\aaa +I(\aaa :\bbb)$ is unmixed follows from Lemma~\ref{unm}.
\end{proof}

\begin{proposition}\label{linkage}
 Let $(R, \m )$ be a local Gorenstein ring and let $I$ be
  an unmixed ideal of height $g < d:= {dim} (R).$ Then for
  every ideal $J$ in the linkage class of $I$ there exist elements $x_1,
  \ldots, x_{d-g}$ such that $(I, x_1, \ldots, x_{d-g})$ is $\m
  $-primary, and for all $t>0$ there exists an ideal $J_t$ in the
  linkage class of $I_t:= (I, x_1^t, \ldots, x_{d-g}^t)$ with
  $J\subset J_t $.
\end{proposition}
\begin{proof}
Let $J$ be an ideal in the linkage class of $I$,
 $$J=\aaa _n : \left(
    \aaa_{n-1} : \ldots :(\aaa _1 :I) \right),$$
with $\aaa _1 \subset I$ and $\aaa _i \subset \aaa _{i-1}:\left(
    \aaa_{i-2} : \ldots :(\aaa _1 :I) \right)$ for $1 < i \le n$, 
Gorenstein ideals of finite projective dimension and of height $g$.
Choose \newline $\bbb \subset \aaa _1 \cap \ldots \cap
\aaa _n$ a parameter ideal of height $g$, and choose  
$\underline{x}=x_1, \ldots, x_{d-g}$ a sys
tem of parameters modulo $\bbb $ (and hence also
a system of parameters modulo each of the $\aaa _i$'s). Then for all
$ t >0$, we claim that
$$J_t:=(\aaa _n , \underline{x}^t):\left( (\aaa _{n-1}, \underline{x}^t): \ldots \left( (\aaa
    _1, \underline{x}^t):(I, \underline{x}^t)\right) \right)
$$
is in the linkage class of ${(I, \underline{x}^t)}$, and 
$J\subset J_t$. Here we use the notation $\und{x}^t=x_1^t, \ldots,
x_{d-g}^t$.

First observe that for every $i=1, \ldots, n $ and for every $t$, the
ideal $(\aaa _i, \underline{x}^t)$ is still a Gorenstein ideal of
finite projective dimension. This is justified by the fact that if
${\bf F.}_i$ denotes a minimal free resolution of $R/\aaa _i$, and ${\bf
K.}(\und{x}^t)$ denotes the Koszul complex of $\underline{x}^t$, then ${\bf
F.}_i \otimes {\bf K.}(\underline x^t)$ is a free resolution of $R/(\aaa _i,
\underline{x}^t)$.
Indeed, it is enough to see that
${\bf
F.}_i \otimes {\bf K.}(\und{x}^t)$ is acyclic. This  is true because its $i$th
homology module is $\mathrm{Tor}_i ^R\left( R/\aaa_i
 , R/(\und{x}^t)\right) $, 
which is zero for all $i\ge 1$ since ${\bf K.}(\underline x^t)$ remains acyclic upon
tensoring with $R/\aaa_i
 $ (the $x$'s form a regular sequence in $R/\aaa _i$).

We prove the claim by induction on $n$. The case $n=1$ is obvious. 
Let $n=2$.

According to Lemma~\ref{mapping cone}, we can write 
$\aaa _2 =\bbb :\delta $
where $\delta \in R$ is such that the last
 map $\psi _g :R \lar R$ in the map of complexes $\psi . $ from the
 Koszul resolution of $R/\bbb $ to a minimal free resolution of
 $R/\aaa _2$ is multiplication by $\delta$.
 Since the corresponding free resolutions for $R/(\bbb, \und{x}^t)
$ and $R/
 ( \aaa _2, \und{x}^t)$
   are obtained from the resolutions of $R/\bbb$, respectively $R/\aaa
   _2$, by tensoring with ${\bf
   K.}(\und{x}^t)$, it follows that the last map of the map of
 complexes between them is again multiplication
   by $\delta $, and therefore we also have $(\aaa _2 , \und{x}^t
)=(\bbb
     , \und{x}^t
):\delta $.
It is enough to show that
$$
(\bbb, \underline{x}^t):\left( (\aaa _1, \underline{x}^t):(I, \underline{x}^t)\right) \supset \bbb :(\aaa
_1 :I),
$$
because then it follows that
$$
(\aaa _2 , \underline{x}^t):[ (\aaa_1, \underline{x}^t):(I,
  \underline{x}^t)] =[(\bbb, \und{x}^t):\delta ]: [(\aaa _1 ,
\und{x}^t):(I, \und{x}^t)]=
$$
$$
[(\bbb , \underline{x}^t):\left( (\aaa _1, \underline{x}^t):(I,
  \underline{x}^t)\right)]:\delta \supset \ \ \ \ \ \ \ \ \ \ \ \ \ \  \ \ \ \
\ \ \ \ \ \ \ \ \ \ \ \  \ \ \ \ 
$$
$$
[\bbb :(\aaa_1 :I)]:\delta= (\bbb :\delta ) :( \aaa_1 :I)=
 \aaa_2 :(\aaa _1:\delta ).\ \ \ \ \ \ \ \ \ \ \ \ \ \ \ 
 \ 
$$
By replacing $\aaa _2$ by $\bbb $ and changing the notation
accordingly, we can assume without loss of generality that $\aaa _2
\subset \aaa _1$. 
By Lemma~\ref{case 1}, in this case we have
$$
\aaa _2:(\aaa_1 :I)=\aaa _2 +I(\aaa _2:\aaa _1)\subset (\aaa _2,
\underline{x}^t)+(I, \underline{x}^t) \left( (\aaa _2,
  \underline{x}^t):(\aaa _1, \underline{x}^t)\right)
$$
$$=
(\aaa _2, \underline{x}^t):\left( (\aaa_1, \underline{x}^t):(I, \underline{x}^t)\right).
$$
This finishes the proof of the claim for the case $n=2$.
Assume $n \ge 3$. Let
$
K:= \aaa _{n-2}: \left( \ldots (\aaa _1:I)\right) 
$
and 
$
K_t:= (\aaa _{n-2}, \underline{x}^t): \left( \ldots \left( (\aaa _1, \underline{x}^t):(I,
    \underline{x}^t)\right) \right). \hfill 
$

The induction hypothesis shows that $K\subset K_t$ for all $t$, hence
\newline 
$(K, \und{x}^t) \subset K_t$.
Therefore we have
$$
J_t=(\aaa _n, \underline{x}^t):[(\aaa_{n-1}, \underline{x}^t):K_t]\supset (\aaa _n,
\underline{x}^t):[(\aaa _{n-1}, \underline{x}^t):(K, \underline{x}^t)],
$$
since double linkage preserves inclusions,
 and by the case $n=2$ applied to $I:=K$ it follows that 
$$J_t \supset \aaa _n :(\aaa_{n-1}:K)=J.
$$
This finishes the proof of the claim and the proof of the proposition.

\end{proof}

We are now ready to prove Thm.~\ref{main}:
\begin{proof}
It is enough to show that $\li{I} (I:\tau )\subset I^*$, since this
will imply that the tight closure of the ideal on the left hand
side is contained in $I^*$, hence
$(\li I)^* (I:\tau )\subset (\li{I} (I:\tau ))^* \subset
I^*.$  

First assume that $I$ is $\m $-primary.

Theorem~\ref{essential} implies that if $d\in R$ is such that
$d^q I^{<q>} \subset I\q $ for all $q$, then $d(I:\tau )\subset I^*$.
We claim that any $d \in \li{I}$ has this property.

It is obvious that any $d \in I$ will have the desired property.
In order to finish
 the proof of the claim, it suffices to prove that
for any 
$\aaa
\subset I$ Gorenstein ideal of finite projective dimension, for all $c \in R$ and for all $q$ we have
$cI^{<q>} \subset I\q \iff cJ^{<q>} \subset J\q $, where $J=\aaa :I$.

Indeed, we have
$$
cI^{<q>} \subset I\q \iff c(\aaa \q :J\q )\subset \aaa \q :J^{<q>}
\iff 
$$
$$
 \aaa \q :J\q \subset \aaa \q:cJ^{<q>} \iff cJ^{<q>} \subset J\q.
$$
Here we used the fact that $I\q = \aaa \q :J\cor $, since $I\q $ is
$\m $-primary.
This shows that for every $K$ in the linkage class of $I$ we have
$cI^{<q>}\subset I\q \iff cK^{<q>} \subset K\q $, and therefore 
for every $d \in K$ we have $d^q I^{<q>} \subset I\q$.

Let $u \in \tilde{I}$ and write
$u=d_1^q +\ldots d_n ^q $ with $d_i \in K_i$, where $K_1,
\ldots K_n$ are ideals in the linkage class of $I$. Since $d_i^q
I\cor \subset I\q$ for all $i$, we have $u^qI\cor \subset I\q$. If $g \in I:\tau
$, it follows that $\tau g^q \in I\cor$, and hence $\tau u^qg^q \in I\q$,
which shows that $ug\in I^*$ as desired.

Now assume that $I$ has height $g <d =\mathrm{dim}(R)$, and let
$J\subset \li{I}$ be
an ideal in the linkage class of $I$. Let $\underline{x}=x_1, \ldots,
x_{d-g}$ and $J_t$ be as in Prop.~\ref{linkage}, for all $t>0$.
According to the $\m $-primary case, we have:
$$
J\left( (I, \underline{x}^t):\tau \right)\subset J_t \left( (I,
  \underline{x}^t):\tau \right) \subset (I, \underline{x}^t)^*
$$
for all $t$. Intersect over all $t$ and use the fact that
$$
\bigcap_t\,  \left( (I, \underline{x}^t):\tau \right) =\left
  ( \bigcap_t\,  (I,
  \underline{x}^t)\right) :\tau =(I:\tau )
$$and
$$
\bigcap _t \, (I, \underline{x}^t)^*=I^*
$$
(the last equality follows because if $u \in \cap_t (I, \und{x}^t)^* $,
then  $cu^q \in \cap _t (I\q , \und{x}^{tq})$ for some
$c \in R^0$ and all $q=p^e$, and so $cu^q \in I\q $ by the Krull
intersection theorem, which means that
$u \in I^*$). It follows that $J(I:\tau )\subset I^*$.

Therefore we have $\tilde{I} (I:\tau ) \subset I^*$, and since the
ideal on the right is tightly closed, it follows that $
  ( \tilde{I})^* (I:\tau )\subset I^*$.
\end{proof}

As an application of this result, we prove an unexpected restriction on the linkage
class of certain ideals:
\begin{corollary}\label{max}
Let $R$ be as above and assume that $I$ is an unmixed ideal containing the
test ideal $\tau$. Then $I^*=(\tilde{I})^*$. 

If in addition $I$ is
tightly closed, 
then $I$ contains every ideal in its linkage class. In particular if
the test ideal $\tau $ is unmixed and tightly closed, then it is
maximal in its linkage class.

If $I$, $J$ are unmixed ideals
containing $\tau $ and are in the same linkage class, then $I^*=J^*$.
\end{corollary}
\begin{proof}
If $I\supset \tau$ then $I:\tau =R$, and
Thm.~\ref{main} implies that $( \tilde{I})^*\subset I^*.$ If in
addition $I$ is
tightly closed it follows that $\li{I} \subset (\li{I})^* =I$, and
therefore $I$ is maximal in its linkage class. The last statement is
immediate, since we have $I^*=J^*=(\tilde{I})^*$.
\end{proof}
The next corollary deals with a property of the test ideal:
\begin{corollary}\label{lit}
Let $(R, \m )$ be a Gorenstein ring of characteristic $p$ with test
ideal $\tau $; assume that
the ideals $\tau $ and $\tau ^2$ are unmixed.
Then we have
$$
\li{\tau } \tau =\tau ^2,
$$
where $\li{\tau }$ is the sum of all ideals in the linkage class of
$\tau $.
\end{corollary}

\begin{proof}
Let $\aaa $ be an ideal generated by  parameters, with
$\aaa \subset \tau ^2$ and of the same height as $\tau $.
Apply the result of Thm.~\ref{main} for $I=\aaa: \tau =\aaa ^*$. Note that
$\li{I}=\li{\tau }$, and $I$ is tightly closed.
We have
$$
\li{\tau } (\aaa :\tau ^2 )\subset \aaa :\tau,
$$
or equivalently
$$
\aaa :\tau ^2 \subset \aaa :\li{\tau } \tau ,
$$
and therefore $\li{\tau }\tau \subset \tau ^2 +\aaa =\tau ^2 $ by the
choice of $\aaa $.
\end{proof}

The final application deals with a particular case of the localization
problem for tight closure:
\begin{corollary}\label{commutes}
Let $R$ be a Gorenstein ring and let $I$ be an unmixed ideal such
$\tilde{I}$ is $\m $-primary.
 Then for any multiplicative system $W$, we have
$(I_W)^*=(I^*)_W$.
\end{corollary}
\begin{proof}
Since $\li{I}$ is $\m $-primary, 
Thm.~\ref{main} shows that $(I:\tau)/I^*$ is a finite length
module and the conclusion follows by Cor.~\ref{testexp}.
\end{proof}
\begin{note}
The hypothesis that $\li{I}$ is $\m-$primary is equivalent to $I_P$ being
in the linkage class of a Gorenstein ideal of finite projective
dimension for every prime ideal
$P \ne \m$ (choose $J_n=\aaa _n : \left(
    \aaa_{n-1} : \ldots :(\aaa _1 :I) \right)$ with $n$ minimal such
that $J_n \not\subset P$. Then
$(J_{n-1})_P =(\aaa_{n-1})_P$, which is in the linkage class of $I_P$,
is a Gorenstein ideal of finite projective dimension). 

With the possible exception of the case when $R$ is an isolated singularity (in which case the localization
of tight closure holds trivially, because the test ideal is $\m
$-primary), it is never the case that the condition
discussed above holds for every ideal in $R$, because it implies
that $I_P$ has finite projective dimension.
Moreover, Thm. 2.10 in ~\cite{PU} gives a large class of non $\m
$-primary ideals that contain every ideal in their linkage classes,
therefore providing examples when the hypothesis of Cor.~\ref{commutes}
fails. Note however that the notion considered in ~\cite{PU} is
linkage by complete intersections; the linkage class by Gorenstein ideals
may be larger.
Thus, the assumption in Cor.~\ref{commutes} is quite restrictive, but this author
believe it is worth adding to the list of cases when tight
closure commutes with localization. 

For a concrete non-trivial example when the hypothesis holds, consider
a codimension two ideal in a regular local ring, such that $R/I$ is
not Cohen-Macaulay, but it becomes Cohen-Macaulay when localized at
any non-maximal prime.
\end{note}

\end{document}